\pdfoutput=1

\documentclass{IEEEtran4PSCC}
\ifCLASSINFOpdf
\else
\fi
\usepackage[cmex10]{amsmath}
\usepackage{mathtools}
\usepackage{hyperref}
\usepackage{graphicx}
\usepackage{amsmath,amssymb}

\makeatletter
\let\old@ps@headings\ps@headings
\let\old@ps@IEEEtitlepagestyle\ps@IEEEtitlepagestyle
\def\psccfooter#1{%
    \def\ps@headings{%
        \old@ps@headings%
        \def\@oddfoot{\strut\hfill#1\hfill\strut}%
        \def\@evenfoot{\strut\hfill#1\hfill\strut}%
    }%
    \def\ps@IEEEtitlepagestyle{%
        \old@ps@IEEEtitlepagestyle%
        \def\@oddfoot{\strut\hfill#1\hfill\strut}%
        \def\@evenfoot{\strut\hfill#1\hfill\strut}%
    }%
    \ps@headings%
}
\makeatother

\psccfooter{%
        \parbox{\textwidth}{\hrulefill \\ \small{21st Power Systems Computation Conference} \hfill \begin{minipage}{0.2\textwidth}\centering \vspace*{4pt} \includegraphics[scale=0.06]{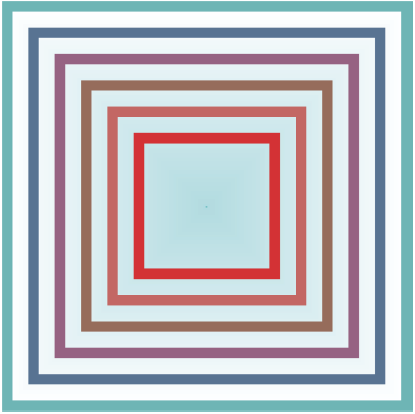}\\\small{PSCC 2020} \end{minipage} \hfill \small{Porto, Portugal --- June 29 -- July 3, 2020}}%
}

\begin{document}
%
\title{Sequential Bayesian Parameter Estimation of Stochastic Dynamic Load Models}

\author{
\IEEEauthorblockN{Daniel Adrian Maldonado, Vishwas Rao, Mihai Anitescu}
\IEEEauthorblockA{Mathematics and Computer Science Division \\
Argonne National Laboratory\\
Lemont, IL, U.S.A\\
\{maldonadod, vhebbur, anitescu\}@anl.gov}
\and
\IEEEauthorblockN{Vivak Patel}
\IEEEauthorblockA{Department of Statistics \\
The University of Wisconsin-Madison\\
Madison, WI, U.S.A\\
vivak.patel@wisc.edu}
}


\maketitle

\begin{abstract}
 In this paper we focus on the parameter estimation of dynamic load models with stochastic terms---in particular, load models where protection settings are uncertain, such as in aggregated air conditioning units. We show how the uncertainty in the aggregated protection characteristics can be formulated as a stochastic differential equation with process noise. We cast the parameter inversion within a Bayesian parameter estimation framework, and we present methods to include process noise. We demonstrate the benefits of considering  stochasticity  in the parameter estimation and the risks of ignoring it. 
\end{abstract}
\begin{IEEEkeywords}
Power System Identification, Power System Dynamics, Load Modeling, Bayesian Statistics
\end{IEEEkeywords}

\thanksto{This material was based upon work supported by the U.S. Department of Energy, Office of Science, under contract DE-AC02-06CH11357.}

\section{Introduction}

Inquiries into fault-induced delayed voltage recovery (FIDVR) events caused by air conditioning units and the impact that a large penetration of user-sized distributed energy resources will have on the dynamic performance of the power grid have prompted much research into developing new load models that can accurately represent the behavior of these devices on the transmission grid. In the United States, these efforts have resulted in the composite load model \cite{NERC2016} and the DER\_A model \cite{dera}.

While these models represent a significant leap in granularity from the previous models used in dynamic stability studies,  development is ongoing. In particular, in order to address the sensitivity issues with block tripping schemes, recent efforts have sought to represent the tripping characteristic of aggregated resources \cite{Ramasubramanian2019, undrill, weber}. New approaches to modeling the tripping characteristics, which often use linear functions of the terminal voltage, pose the question of how to estimate the parameters that define these functions.

Parameter estimation of dynamic load models is still an open question; and since validation remains challenging, disparate approaches coexist. On the one hand, the North American Electric Reliability Corporation (NERC) is actively working on extrapolating data characterized by regions and climate zones to populate the composite load model \cite{NERC2017}. On the other hand, the research impetus in machine learning and deep learning techniques has permeated into the field of power systems, promising to extract insights from the increasing amount of data \cite{zhang2017}. Moreover, the increase in computational power has opened the posibility of elaborating dynamic load models via co-simulation \cite{Chaspierre2018}. These approaches offer valuable contributions and indicate that revisiting the foundations of the load modeling practice is a worthwhile endeavor.

In this paper we revisit the topic of dynamic load parameter estimation with the additional complexity of progressive tripping. Inspired by recent results that use stochastic models to represent short-term load behavior \cite{Milano2013, Roberts2016}, we propose a model that includes uncertainty in the tripping process, resulting in a stochastic dynamic load model. Then, building on previous work \cite{Maldonado2017}, we introduce techniques to perform Bayesian parameter estimation in stochastic models. In Section \ref{section: tripping} we introduce the use of process noise to reflect the uncertainty of the aggregated tripping mechanism. In Section \ref{section: paramest} we pose the mathematical formulation of the parameter estimation problem with system uncertainty or process noise. In Section \ref{section: bayesian} we introduce the Bayesian estimation methodology for load models with system uncertainty. In Section \ref{section: case} we present a case study of the model in both deterministic and stochastic form. In Section \ref{section: conc}  we summarize our conclusions.

%

\section{Load Progressive Tripping Models and Stochastic Tripping Characteristics}
\label{section: tripping}

Recent research has pointed out the importance of correct modeling of the protection settings in load models to better understand their impact in stability simulations. Two  salient cases are the study of the protection settings in behind-the-meter photovoltaic (PV) panels and in induction motors to understand fault-induced delayed voltage recovery events. Protection action has been traditionally modeled as ``block tripping,'' and it can be represented with the following equations:
\begin{align}
L_{\textit{frac}} &= \begin{cases}
1 &\text{if } V \geq V_{\textit{1off}} ,\\
0 &\text{if } V \leq V_{\textit{2off}} ,\\
B_{\textit{frac}} &\text{if } V \in (V_{\textit{1off}}, V_{\textit{2off}}) ,\\
\end{cases}
\end{align}
where a fixed percentage $B_{\textit{frac}}$ of the load is tripped when the terminal voltage drops below a certain threshold. Some researchers have argued, however, that this model does not represent the true behavior. Whereas a voltage drop to $V_{\textit{1off}}$ would result in no load loss, a voltage drop to $V_{\textit{1off}} - \epsilon$ would result in a significant tripped fraction. Since the dynamic load models  represent the aggregated action of many devices, subjected to different terminal voltages, that might not trip at the same time, this model has been deemed unrealistic.

To overcome the sensitivity problem, researchers have proposed a series of \textit{progressive tripping} models  to make the tripping fraction of the load a smoother function of the voltage. In its simplest form, a linear characteristic that emanates from the feeder topology is used to compute a tripped load fraction that then is passed through a lag block.
\begin{align}
d_{\textit{input}} &= \begin{cases}
1 &\text{if } V \geq V_{\textit{1off}} ,\\
0 &\text{if } V \leq V_{\textit{2off}} ,\\
\frac{V - V_{\textit{2off}}}{V_{\textit{1off}} - V_{\textit{2off}}} &\text{if } V \in (V_{\textit{1off}}, V_{\textit{2off}}) ,\\
\end{cases} \\
fr_{k + 1} &= fr_k + \frac{h}{T_d}(-fr_k + d_{\textit{input}}) \label{eq:lag_det}
\end{align}

Unfortunately, although approximating the progressive tripping characteristic with a linear function produces qualitatively reasonable results, it fails to represent the true tripping characteristic \cite{Takenobu2018}. One could introduce more complex characteristics such as higher order polynomials, but to obtain parameters for such equations would be difficult. Since the distribution network is a complex system and the disconnected fraction is a parameter that evolves conditioned to many parameters that vary over time, another approach is to consider part of the disconnection behavior as uncertain and represent it with a stochastic process.

\begin{figure}[h]
\centering
\includegraphics[width=0.45\textwidth]{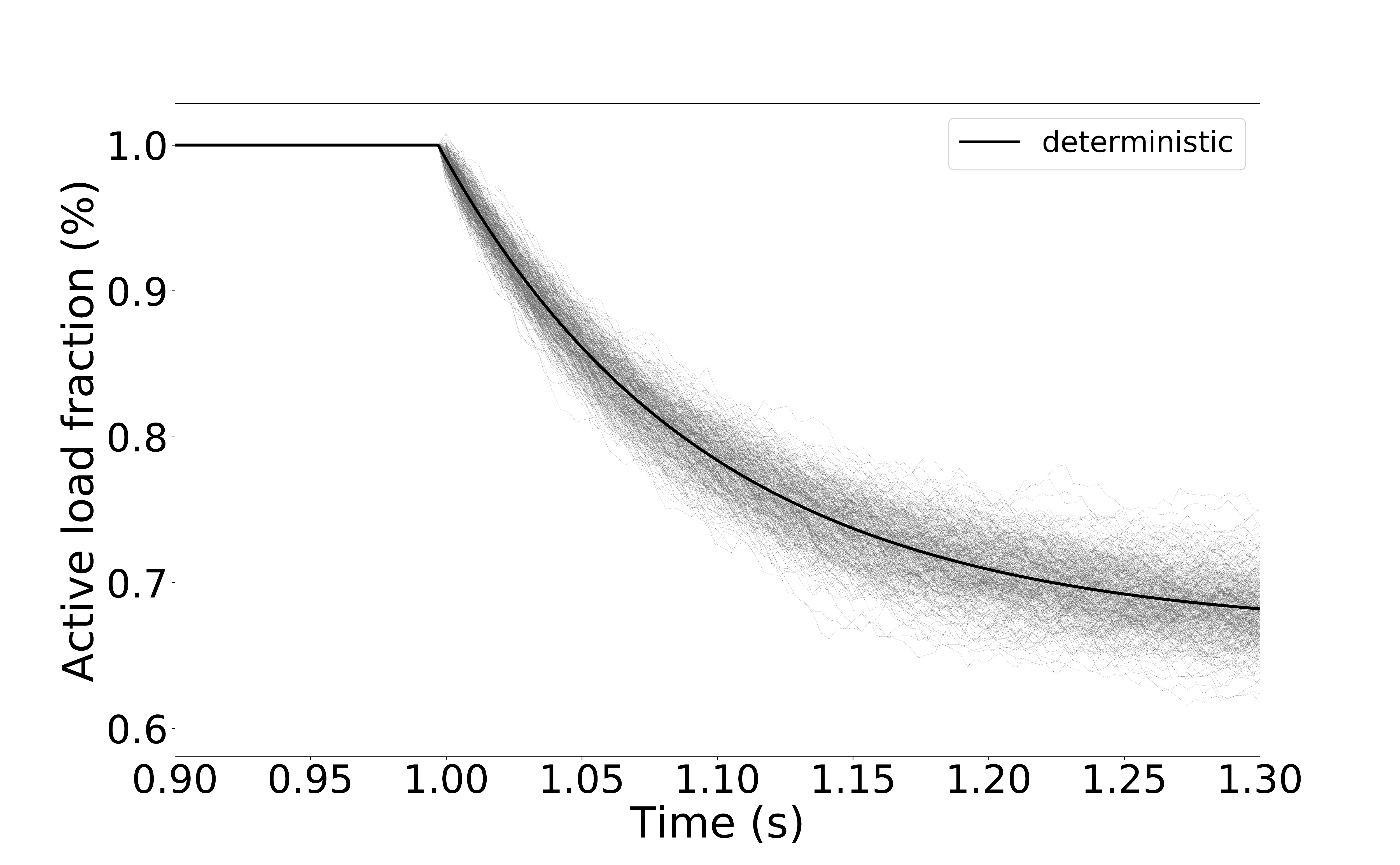}
\caption{Stochastic tripping.}
\label{fig:strip}
\end{figure}

For $V \in (V_{\textit{1off}}, V_{\textit{2off}})$ we can write Equation (\ref{eq:lag_det}) with a stochastic Wiener process.
\begin{align}
fr_{k + 1} &= fr_k + \frac{h}{T_d}(-fr_k + d_{\textit{input}}) +  \Delta W \,,   
\end{align}
where $\Delta W = W_{k+1} - W_{k}$ is $\mathcal{N}(0, \Delta)$ increment of the Wiener process $W$.

The effects  can be seen in Figure \ref{fig:strip}. The stochastic process essentially represents \textit{the part of the model behavior that we cannot explain}. This term is common in most filtering and estimation literature and is often referred to as model error or system noise.

\section{Parameter Estimation}
\label{section: paramest}

The  nonlinear stochastic model can be expressed within the standard Bayesian dynamical model framework:
\begin{align}
x_k &= f(x_{k-1}, \theta, v_k) \label{eq:propagation} \,, \\
y_k &= g(x_k, \theta, w_k) \label{eq:measurement} \,,
\end{align}
with unknown fixed parameters $\theta$. In this model, $x_k \in \mathbb{R}^n$ is the hidden state of the system at time $k$, $y_k \in \mathbb{R}^m$ the measurement, $v_k$ the process noise, and $w_k$ the measurement noise. This model can be also interpreted from a probabilistic perspective  \cite{Cappe2007}, where (\ref{eq:propagation}) can be written as the \textit{transition density}, $p(x_k|x_{k-1};\theta)$ and (\ref{eq:measurement}) as the \textit{observation density}, $p(y_k|x_k;\theta)$. These equations implicitly carry the Markovian assumption. 

The parameter that satisfies the set of state and measurement equations \eqref{eq:propagation} and \eqref{eq:measurement} can be obtained by multiple methods. One can use a variational approach problem within an optimization framework and obtain the parameter $\theta$ that best fits the observations while satisfying the dynamics. This solution is also referred to as maximum a posteriori  (MAP) estimation. Since both the observations and the process is noisy, a single value of $\theta$ is insufficient to match different realizations of the observations. Instead, one looks for a distribution, referred to as posterior distribution, that can suitably explain multiple realizations of the noise. For special problems, one can also obtain the associated uncertainty in the MAP estimate in the variational approach. For example, when $f$ and $g$ are linear, the observation noise is Gaussian, and there is no process noise, then the inverse of the Hessian at the MAP point is also the covariance of the associated MAP. However, $f$ and $g$ are rarely linear, and process noise is common in most scenarios. To address this shortcoming of the variational approach, we use a fully Bayesian approach to mitigate the effects of nonlinearity and the presence of process noise. The Bayesian approach in general attempts to describe the posterior distribution. Typically, unnormalized density of the posterior distribution can be evaluated, and in such scenarios the standard approach is to draw samples from the posterior distribution by using a Markov chain Monte Carlo method (MCMC). These samples help us characterize the posterior distribution or the uncertainty associated with parameter estimates.

\section{Bayesian Parameter Estimation}
\label{section: bayesian}
Bayesian estimation combines the predictions from an analytic model of the system with measurement data, with the goal of estimating a model parameter $\theta$. The core of the Bayesian framework is that it presupposes some \textit{prior} knowledge about the parameter $\theta$, which we call prior distribution $p(\theta)$, and combines it with data $y_{1:N}$ to obtain \textit{a posteriori} knowledge or posterior distribution $p(\theta \mid y_{1:N})$. Furthermore, the inference process includes a ``best" estimate of the parameter along with associated uncertainties.
Equations (\eqref{eq:propagation}) and (\eqref{eq:measurement}) can be written as transition densities:
\begin{align}
  x_N \mid x_{N-1} &\sim f_{\theta}(x_N \mid x_{N-1})\,, \label{eq:transition_dynamics}\\
  y_N \mid x_N &\sim g_{\theta}(y_N \mid x_{N})\,, \label{eq:transition_measurement}
\end{align}
where $f_{\theta}$ and $g_{\theta}$ denote the transition densities for the state and observations respectively, for a static parameter $\theta$. The initial distribution of the process is characterized by $x_0 \sim \mu_{\theta}(\cdot)$. A fully Bayesian paradigm requires specification of the likelihood $p_(y_{1:N} \mid \theta)$ and a prior $p(\theta)$. The unnormalized posterior density is given by the product of the likelihood and prior distributions:
\begin{align}
p(\theta \mid y_{1:N}) &\propto  p(y_{1:N} \mid \theta) \times p(\theta) \label{eq:posterior}\,.
\end{align}
Thus, by characterizing the posterior distribution $p(\theta \mid y_{1:N})$ we can perform statistical analyses (mean, variance, etc.) that provide information about our certainty of the parameter $\theta$. In many cases such as ours, the posterior $p(\theta \mid y_{1:N})$ does not have a closed-form expression, and we have to resort to sampling methods such as MCMC.

\subsection{Markov Chain Monte Carlo Approach}\label{sec:MCMC}
MCMC belongs to a class of methods where a Markov chain is built in such a way that the equilibrium distribution of the Markov chain is the same as the desired posterior distribution. Once the posterior distribution is obtained, various integrals associated with the posterior distribution (such as expectations and covariance) can be computed by using the Monte Carlo integration technique. A Markov chain is a sequence of random variables $X_1, X_2, \cdots, X_n$ such that the conditional distribution of $X_{n+1}$ depends only on $X_n$, which  can be written mathematically  as
\begin{align}
  &\textrm{Prob}(X_{n+1} = \xi \mid X_n =\xi_n, X_{n-1} = \xi_{n-1}, \cdots, X_1 = \xi_1) = \nonumber\\
   &\textrm{Prob}(X_{n+1} = \xi \mid X_n =\xi_n) \,.
\end{align}
Metropolis-Hastings (MH) is the most popular MCMC algorithm \cite{Metropolis1953,Hastings1970}. To sample from a distribution $\varphi(\xi)$, the MH algorithm constructs a transition kernel to go from state $\xi_i$ to $\xi_j$ by a two-step process:
(i) specify a proposal distribution $q(\xi_j \mid \xi_i)$ and (ii) accept draws from  $q(\xi_j \mid \xi_i)$ with an acceptance ratio $\alpha(\xi_i, \xi_j) = \textrm{min} \left \lbrack 1, \frac{\varphi(\xi_j) q(\xi_i \mid \xi_j)}{\varphi(\xi_i)q(\xi_j \mid \xi_i)}\right \rbrack $. 

For our specific problem, we need to sample from $p(\theta \mid y_{1:N})$. To do so  requires that we evaluate  $p_{\theta}(y_{1:N})$ for a proposed $\theta$; and based on MH acceptance criterion, the proposed $\theta$ is accepted or rejected. This process has to be repeated for multiple values of $\theta$ in order to generate a Markov chain. Thus, the likelihood can be rewritten as
\begin{align}
  p(y_{1:N} | \theta) =  p_{\theta}(y_{1:N} \mid x_{1:N}) \times  p_{\theta}(x_{1:N})\,. \label{eq:likelihood_conditional}
\end{align}

The terms $p_{\theta}(x_{1:N})$ and $ p_{\theta}(y_{1:N} \mid x_{1:N})$ can be evaluated by using \eqref{eq:transition_dynamics} and \eqref{eq:transition_measurement}, respectively. However, we have process noise, as described in \cite{SvenssonSL:2018}, and therefore need to integrate over all possible trajectories:

\begin{align}
  p(y_{1:N} | \theta) =  \int_{X^{T+1}} p_{\theta}(y_{1:N} \mid x_{1:N})  p_{\theta}(x_{1:N}) dx_{0:T}\,. \label{eq:likelihood_conditional_stochastic}
\end{align}

To evaluate this integral, we need to sample the process dynamics. In other words, we need to simulate trajectories of the model and evaluate the likelihood until  the integral converges. The convergence  can be slow, especially as the time horizon increases. Therefore, several strategies have been developed to make the process computationally effective.

\subsubsection{Monte Carlo evaluation} \label{sssec:mc}

A simple approach for approximating $p_{\theta}(y_{1:N})$ is to use a Monte Carlo approximation by averaging over different trajectories of $x$. In practice, this involves generating $l$ trajectory samples by integrating the model equations $L$ times (each one will be distinct because of the process noise) and evaluating

\begin{align}
  p_{\theta}(y_{1:N}) \approx \frac{1}{L} \sum_{i=1}^{L} p_{\theta}(x^i_{1:N}) \times  p_{\theta}(y_{1:N} \mid x^i_{1:N}) \,. \label{eq:mclieklihood}
\end{align}

This method deteriorates as the time horizon increases as it becomes harder to explore the proposal distribution. To overcome this drawback, researchers have proposed \cite{Andrieu2010,SvenssonSL:2018} the use of sequential Monte Carlo (particle filter) to sample more efficiently from the proposal distribution.

\subsubsection{Particle Filter evaluation} \label{sssec:pf}

Sequential Monte Carlo, or particle filter, exploits the temporal structure of the probability distribution. In this approach one generates samples $\{x^i_k\}^{\ell}_{i=1}$ for $k = 1, \cdots, N$ such that the samples with index $k$ are approximately distributed according to $p_{\theta}(x_k \mid y_{1:k-1})$. We note that in the previously described Monte Carlo approach, $x_{1:N}^i$ is drawn independent of the measurements $y$. In the particle filter approach, however, the samples interact between time steps. The empirical distribution with $\ell$ samples $\{x^i_k\}^{\ell}_{i=1}$ that approximates $p_{\theta}(x_k \mid y_{1:k-1})$ can be written as
\begin{align}
  p^{PF}_{\theta}(x_k \mid y_{1:k-1}) = \frac{1}{N} \sum_{i=1}^{\ell} \delta_{x^i_k}(x_k)\,,
\end{align}
where $\delta_{x^i_k}(x_k)$ is a point-mass distribution at $x^i_k$. The samples at time step $k=0$ are obtained from the prior, and the samples from $p_{\theta}(x_{k+1} \mid y_{1:k})$ are obtained by propagating the existing samples from  $p_{\theta}(x_{k} \mid y_{1:k-1})$ using the dynamics (Equation \eqref{eq:propagation}). Extensive details about the particle filter approach can be found in \cite{SvenssonSL:2018,Andrieu2010}.

\section{Case Study}
\label{section: case}
For simplicity in both implementation and exposition, we use the model described in \cite{Renmu2006} consisting of a ZIP model together with a third-order induction motor, which can be considered  a subset of the WECC load model. Our modification to the equation is that the current injection of the motor is multiplied by the fraction of load computed by the progressive tripping scheme. The load voltage is subjected to a transient following the voltage test function described in \cite{dera}:
\begin{align}
\label{eq:voltage}
V(t) &= \begin{cases}
a &\text{if } 1 \leq t < (1 + \frac{b}{60}) ,\\
\frac{-(1-d)}{\frac{b}{60} - c}(t - (1 + c)) + 1 &\text{if } (1 + \frac{b}{60}) \leq t < 1 + c ,\\
 1 &\text{otherwise}
\end{cases}
\end{align}  
In this case study we are concerned with estimating the voltage thresholds of the progressive tripping schemes, together with the inertia of the motor load. We  first tackle the case in which the underlying model is deterministic; that is, we assume the model that we have is a perfect representation of the underlying system except for our ignorance about the parameters. We perform the Bayesian parameter estimation with MCMC, and we test the results under different voltage depressions. Then, we move to the stochastic case, in which we assume we cannot explain part of the behavior of the model, in this case the tripping mechanism. We compare what effect the uncertainty of the load model has on the parameter estimation results. We show the capability of the particle filter MCMC to compute the results efficiently. For the MCMC sampling we use the open source library \textit{emcee} \cite{ForemanMackey2013}. We implement the likelihood computation via MC and stochastic MC with our own code, where the equations in the Appendix A are integrated with a forward Euler scheme. We note that in the MCMC framework we work with the logarithm of the likelihood, or log-likelihood.

\subsection{Deterministic Model}
\begin{figure}[h]
\centering
\includegraphics[width=0.45\textwidth]{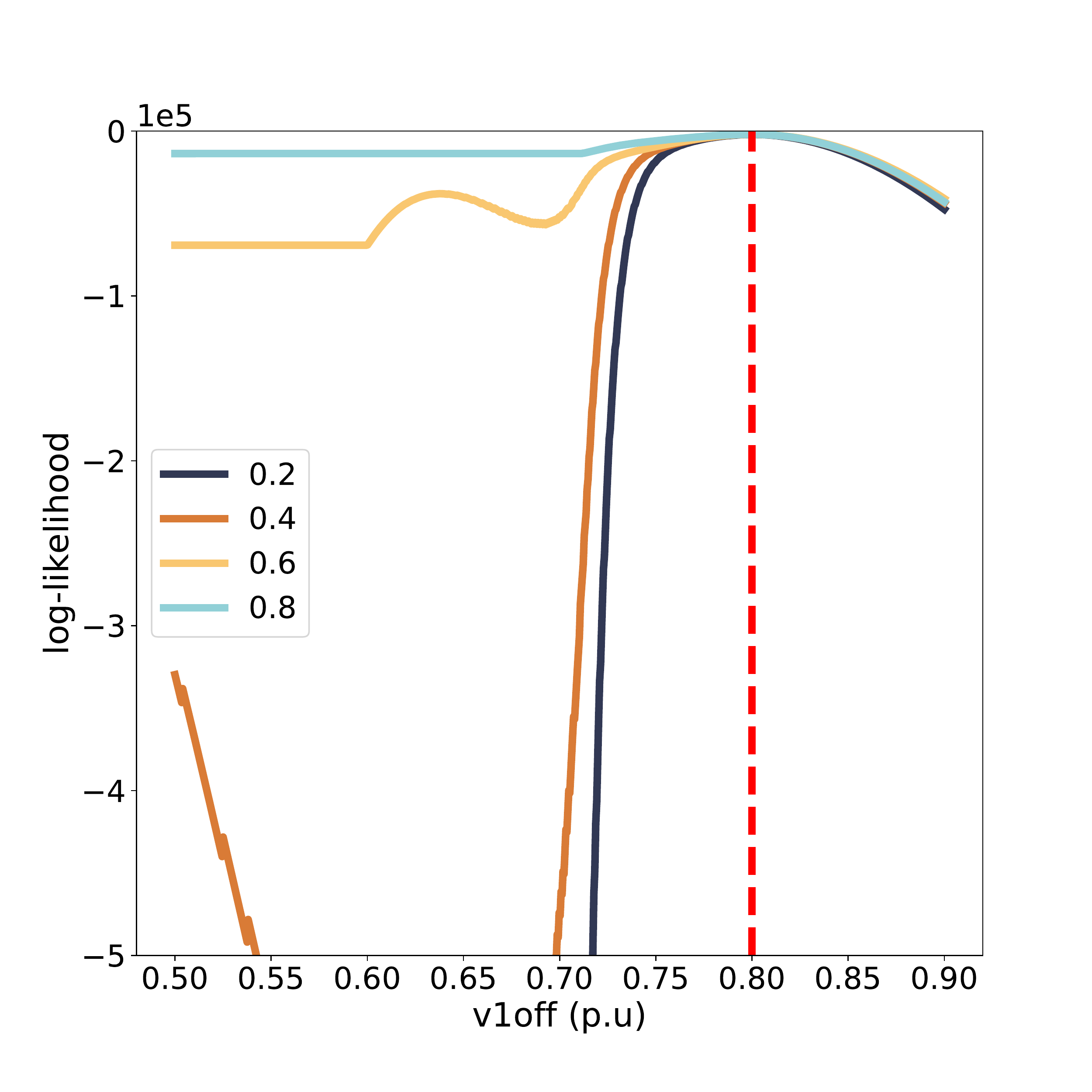}
\caption{Log-likelihood $l(y_{1:t} \mid V_{1\textit{off}})$ of the deterministic system for different voltage drops. Red represents the true value.}
\label{fig:detlike_v1off}
\end{figure}

\begin{figure}[h]
\centering
\includegraphics[width=0.45\textwidth]{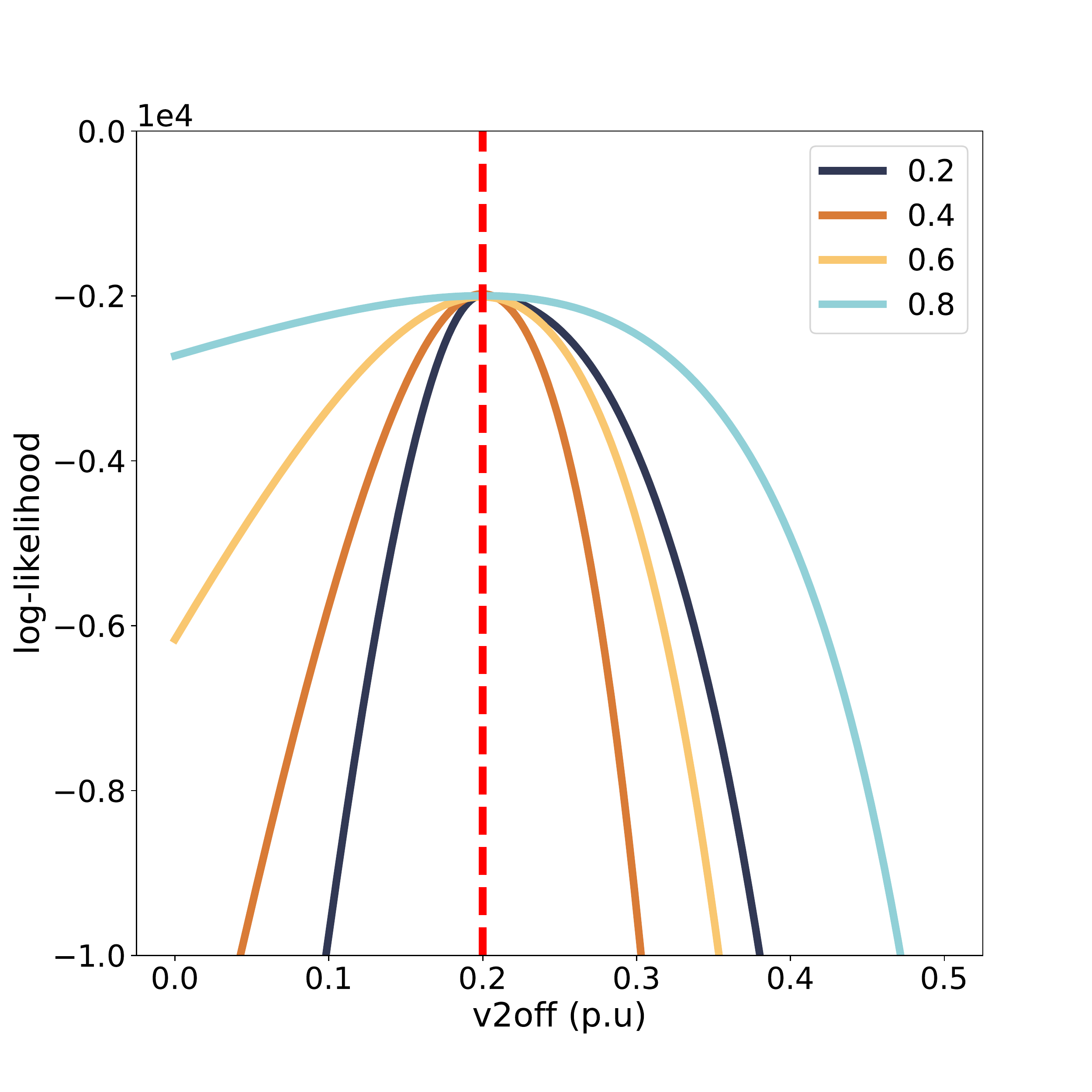}
\caption{Log-likelihood $l(y_{1:t} \mid V_{2\textit{off}})$ of the deterministic system for different voltage drops. Red represents the true value.}
\label{fig:detlike_v2off}
\end{figure}

For the first experiment we consider the underlying system to be deterministic. Here, we consider as unknowns the parameters $V_{1\textit{off}}$, $V_{2\textit{off}}$, and $H$. We use (\ref{eq:voltage}) to generate a voltage drop to $0.6$ pu. by setting the parameter $a$ to this value. The measurement data is generated by applying this voltage to the load model described in  Appendix \ref{sec:load model} and adding Gaussian noise of variance $0.01$ to the active and reactive power measurements. 
Before running the MCMC sampling, we visualize the shape of the log-likelihood functions as we vary individual parameters.  In Figures \ref{fig:detlike_v1off} and \ref{fig:detlike_v2off} we plot the value of the log-likelihood for different values of $V_{1\textit{off}}$ and $V_{2\textit{off}}$,  keeping the rest of the parameters to the true value.  

\begin{figure}[h]
\centering
\includegraphics[width=0.45\textwidth]{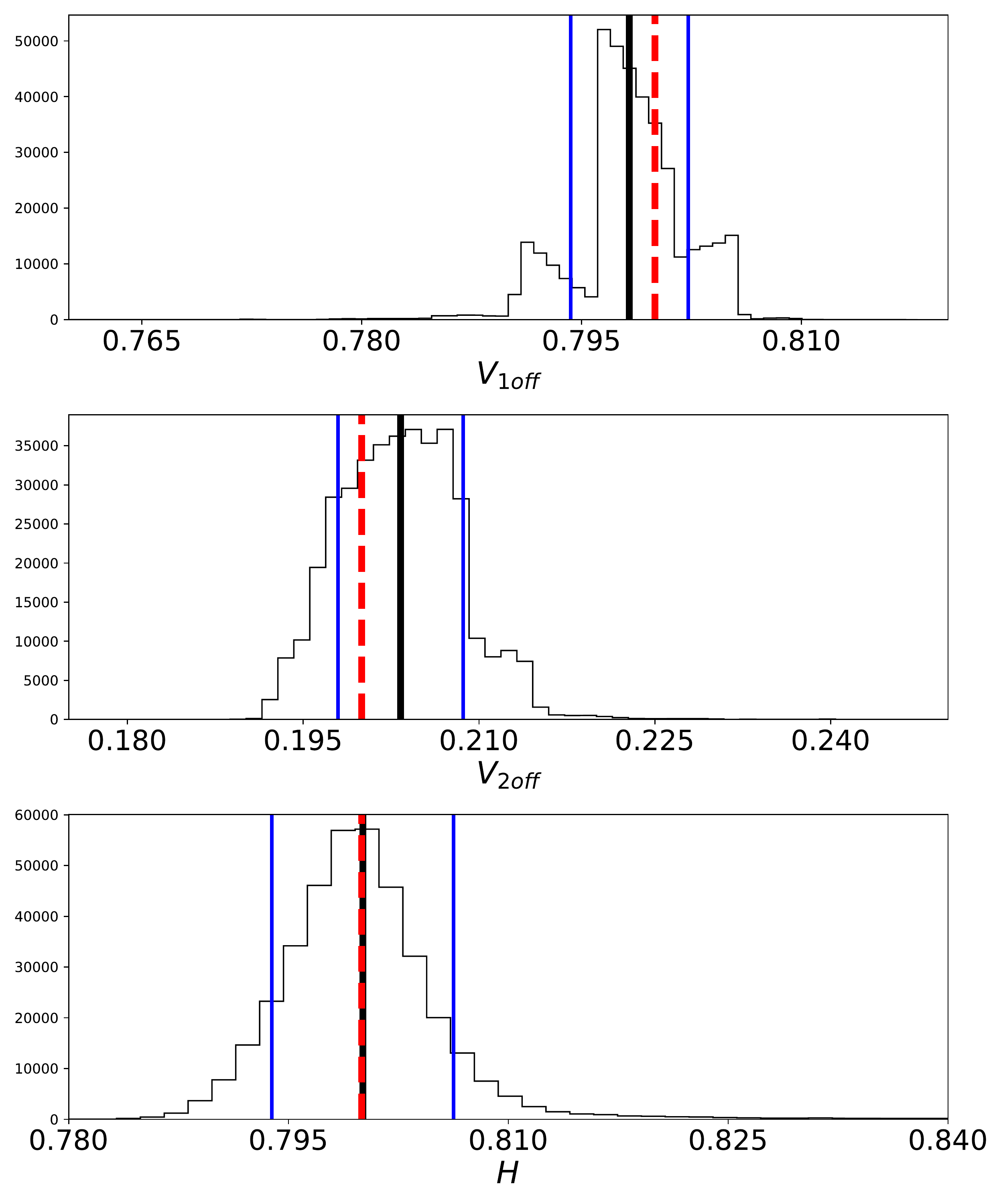}
\caption{Posterior distributions for the deterministic case. The red vertical line represents the true value. The black vertical line represents the mean of the posterior. Blue vertical lines represent one standard deviation of the posterior.}
\label{fig:posterior_deterministic}
\end{figure}

We can see that as the voltage minimum decreases, the log-likelihood values around the true value decrease, making them less likely and pointing to the true value. This situation is to be expected: more acute voltage drops will result in richer data that allow us to determine the slope of the tripping characteristic. We also note the presence of various local minima in $l(y_{1:t} \mid V_{1\textit{off}})$ and its strong non-linearity. We perform the MCMC sampling by initializing the chain around biased values. We choose flat priors for $V_{1\textit{off}}$, $V_{2\textit{off}}$, and $H$ that range from $0.5$ to $0.9$, from $0.1$ to $0.3$, and from $0.7$ to $1.1$, respectively. The chain consists of $400$ walkers and $1,000$ steps. Figure \ref{fig:posterior_deterministic} shows representations of the posterior distribution by plotting histograms of the MCMC chains. We can see that whereas the posterior distribution of the motor inertia $H$ seems Gaussian-like, the posterior distributions for the tripping characteristic parameters present distributions that would be difficult to represent parametrically. Regardless, the posterior distributions obtained with MCMC seem to represent the true values satisfactorily.

\subsection{Stochastic Model}

\begin{figure}[h]
\centering
\includegraphics[width=0.5\textwidth]{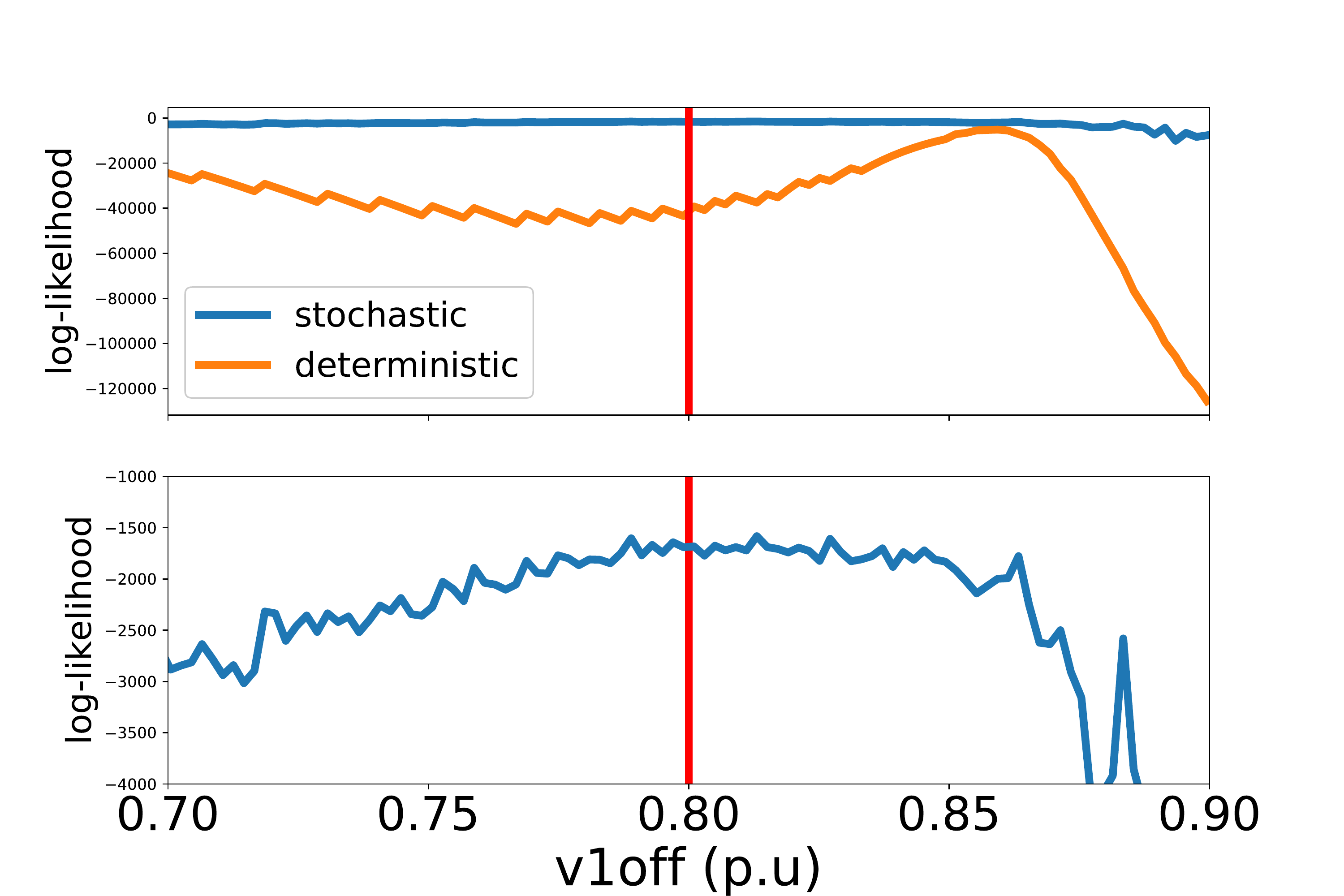}
\caption{Log-likelihood $p(y_{1:t})$ considering deterministic system and stochastic system. Because of the system error, there is a bias in the deterministic prediction.}
\label{fig:loglike_stoch}
\end{figure}

For the stochastic model we perform the same experiment as before but now the load model is integrated with the addition of process noise. The process noise is Gaussian (white) with variance $0.01$. In Figure \ref{fig:loglike_stoch}, as in the preceding subsection, we show a plot of the log-likelihood, but this time we plot the log-likelihood computed without system noise (\ref{eq:likelihood_conditional}) and integrating over the system noise (\ref{eq:likelihood_conditional_stochastic}) for the parameter $v_{1\textit{off}}$. We call these deterministic log-likelihood and stochastic log-likelihood, respectively. Note that the maximum of the deterministic log-likelihood no longer coincides  with the true value because the process noise introduces bias. The maximum of the stochastic likelihood still seems to match the true value. However, we can see that the area around the true value seems to be  flatter compared with the deterministic case, implying that they are as likely as the true value. Intuitively, we might say that when we average over the process noise, small deviations around the tripping parameters do not matter that much.

\begin{figure}[h]
\centering
\includegraphics[width=0.45\textwidth]{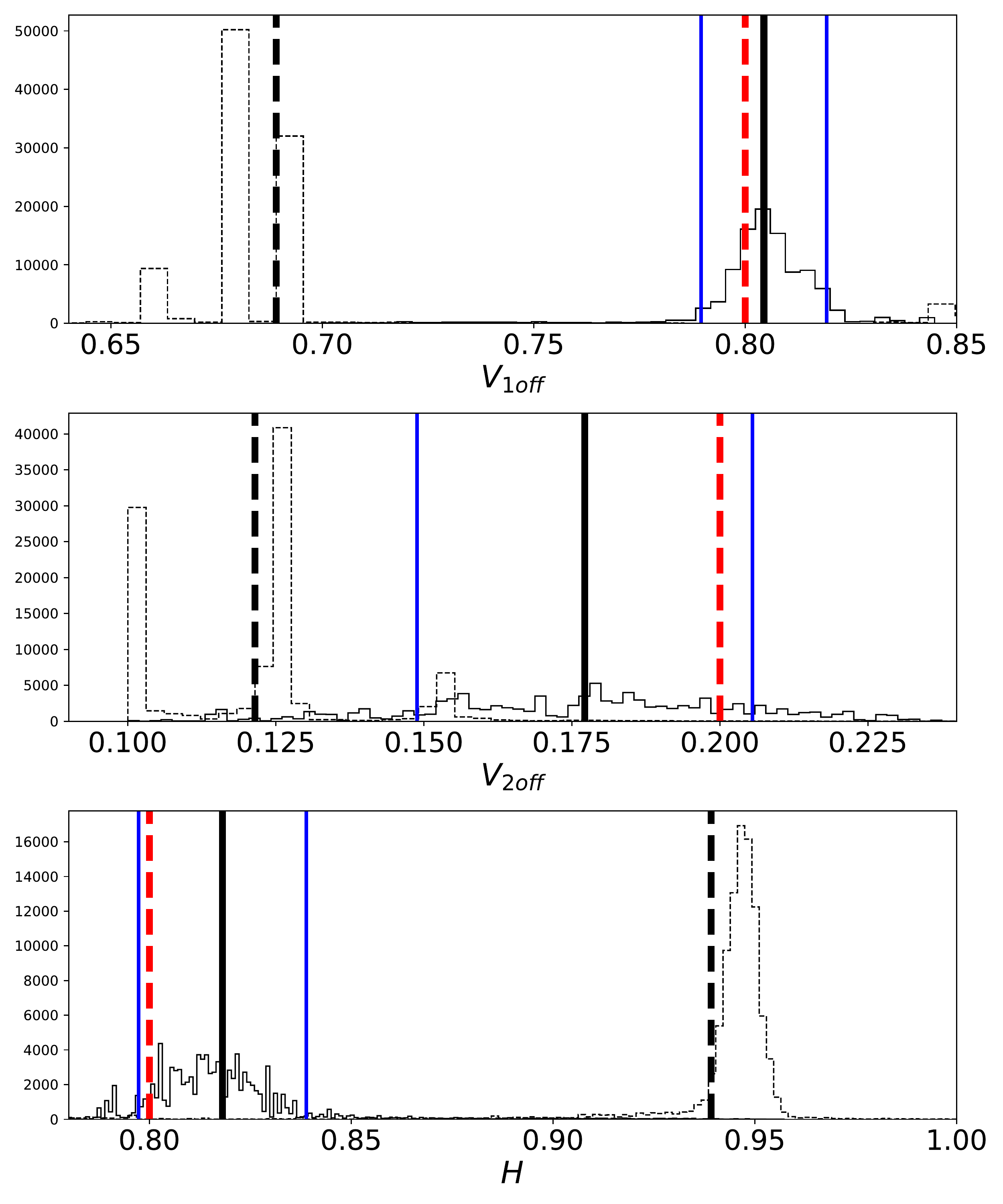}
\caption{Deterministic posterior distribution (discontinuous histogram) with mean (discontinuous black vertical line), Stochastic posterior distribution (continuous histogram) with mean (continuous black vertical line) and standard deviation (blue vertical lines). Red represents the true value.}
\label{fig:posterior_stochastic}
\end{figure}

\begin{figure}[h]
\centering
\includegraphics[width=0.5\textwidth]{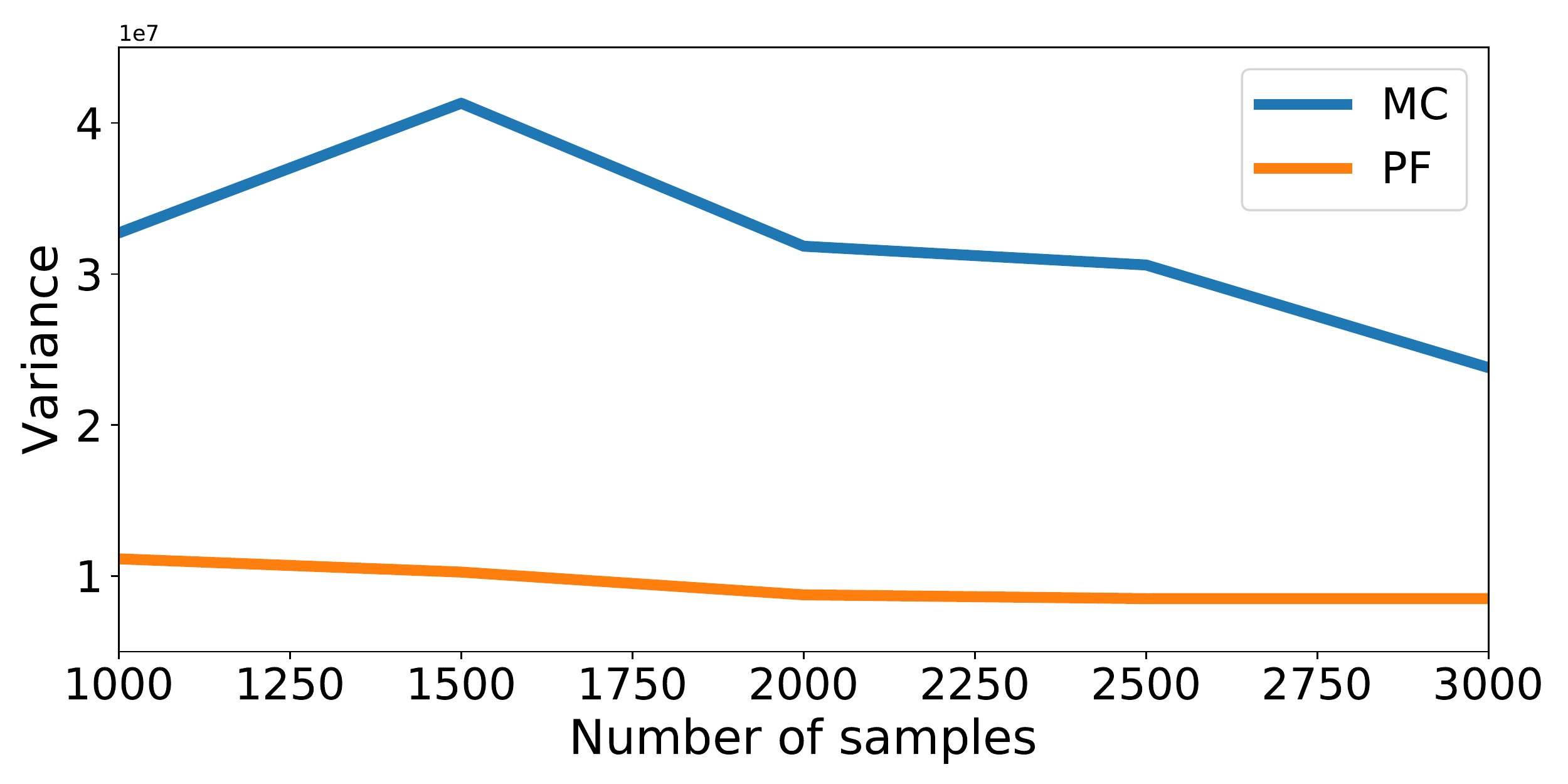}
\caption{In this plot, the faster convergence of the particle filter is shown. We plot the variance of the log-likelihood function over a range of values, as we increase the number of samples. We can see how the variance of the Monte Carlo approach decreases more slowly.}
\label{fig:pfconv}
\end{figure}

In Fig.~\ref{fig:posterior_stochastic} we show again the results of MCMC sampling. This time we draw samples from both the deterministic posterior (neglecting process error) and the stochastic posterior (integrating over all trajectories). The mean of the posterior that has been sampled by using a deterministic likelihood presents an important bias with respect to the true value. The mean of the posterior that has been sampled by integrating the process noise, however, has a mean value closer to the true value, but  the posterior distribution presents much higher variance. Hence, we can say that introducing process noise increases the variance of our estimates. Finaly, in Fig.~\ref{fig:pfconv} we show the different convergence rates of the log-likelihood computed with Monte Carlo (Section \ref{sssec:mc}) and Particle Filtering (Section \ref{sssec:pf}) by plotting the variance of the log-likelihood function over a range of values as we increase the number of samples (particles). Whereas the case that we are tackling is relatively simple and low-dimensional (we only consider three parameters, $v_{1\textit{off}}$, $v_{2\textit{off}}$ and $H$), as we increase the number of parameters, the computations become much more onerous and it is necessary to resort to methods such as particle filtering.

\section{Conclusions}
\label{section: conc}
In this work we have discussed recent developments in aggregated dynamic load modeling, and we have proposed including a stochastic term in the tripping characteristic. Whereas many distinct models  produce satisfactory results with regard to reproducing events, we believe that the main utility of our model is its  power to explore what part of the system response can be explained by the model. Being able to model the uncertainty and include this in the estimation process is paramount to obtaining more robust and generalizable estimates.

Previously, when we estimate a number of parameters, we have assumed that the rest of the parameters are known. This assumption, although not realistic, was done for simplification, but it merits a short discussion. Estimating the parameters of the complex load model can be abstracted as follows: Given a function from a high-dimensional space to a low-dimensional space:
\begin{equation}
x = f(\theta_1, \theta_2, \dots, \theta_n) \,,
\end{equation}
where $x$ is a scalar, obtain $\theta_i$, given evaluation points $x_j$. In general, many distinct combinations of $\theta_i$ might minimize this functional relationship, which means that the problem can be considered as  being \textit{ill-posed}. Furthermore, if the points $x_i$ are confined to a region, the $\theta_i$'s that minimize this might not be the same as the ones that minimize the observations over a larger region. This issue, overfitting, is a common problem with high-dimensional models. As we increase the complexity of models with the aim of better capturing the underlying physics of the load, we must also consider what  the added complexity will entail for parameter estimation.

Complex load models are a step forward because they attempt to develop the behavior of the load from first principles. One cannot, however, invert their parameters simply from measurements. The initiative by NERC on extrapolating parameters by regions is a necessary step to populate this models. However, we must not reject the wealth of measurements that diverse, high-fidelity sources such as $PMUs$ give us. Bayesian statistics can provide a scientifically robust framework where extrapolated parameters can be set as prior distributions and be combined with high-frequency measurements of events.

\begin{appendices}

\section{Load model}
\label{sec:load model}

For load representation it is standard to use the ZIP load model:
\begin{align}
P_{\textit{zip}}(V) &= P_p + P_i \left( \frac{V}{V_0}\right)  + P_z  \left(\frac{V}{V_0} \right)^2 \,, \\
Q_{\textit{zip}}(V) &= Q_p + Q_i \left( \frac{V}{V_0}\right)  + Q_z \left(\frac{V}{V_0} \right)^2 \,,
\end{align}
where $P_{\textit{ZIP}}$ and $Q_{\textit{ZIP}}$ are the total ZIP demand. The motor equations are
\begin{subequations}
\begin{align}
\dot{e'_{d}} &= \frac{-1}{T_p}(e'_d + (x_0 - x')i_q) + s \omega_s e'_q \,, \\
\dot{e'_{q}} &= \frac{-1}{T_p}(e'_q - (x_0 - x')i_d) - s \omega_s e'_d \,, \\
\dot{s} &= \frac{1}{2H}(\tau_m - e'_d i_d - e'_q i_q) \,, \\
0 &= r_a i_d - x' i_q + e'_d + V sin(\theta) \,, \\
0 &= r_a i_q - x' i_d + e'_q  - V cos(\theta) \,.
\end{align}
\end{subequations}

The active power and reactive power  consumed by the motor are written respectively  as
\begin{subequations}\begin{align}
P_{\textit{mot}} &= - V sin(\theta) i_d + V cos(\theta) i_q \,, \\
Q_{\textit{mot}} &=   V cos(\theta) i_d + V sin(\theta) i_q \,.
\end{align}\end{subequations}

The measurements are the result of the ZIP load and the power consumed by the motor, adjusted by the active fraction from the tripping characteristic of Equation (\ref{eq:lag_det}):
\begin{align}
P_{\textit{inj}}(V, t) &= frP_{\textit{mot}}(V, t) + P_{zip}(V) \,, \\
Q_{\textit{inj}}(V, t) &= frQ_{\textit{mot}}(V, t) + Q_{zip}(V) \,.
\end{align}

\begin{figure}[h]
\centering
\includegraphics[width=0.45\textwidth]{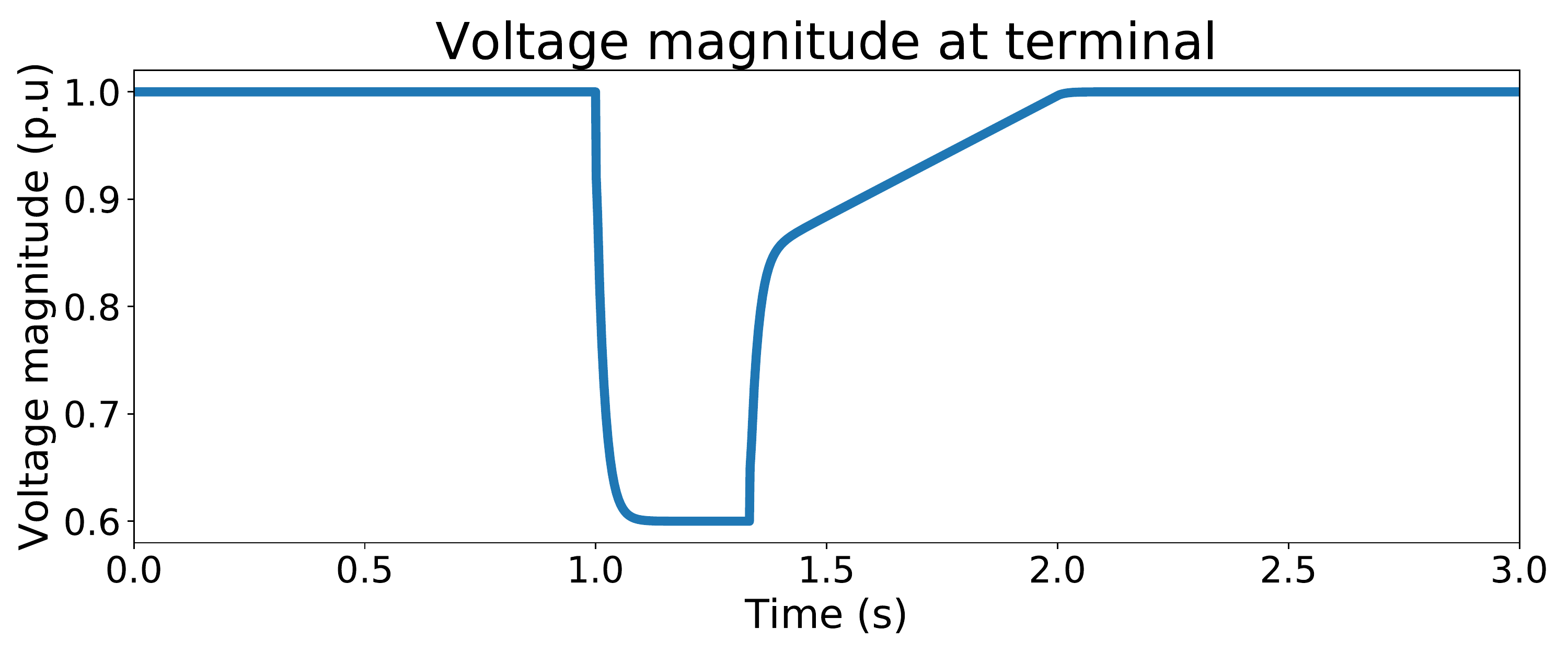}
\caption{Voltage magnitude disturbance at load terminal.}
\label{fig:voltage_dist}
\end{figure}

\begin{figure}[h]
\centering
\includegraphics[width=0.48\textwidth]{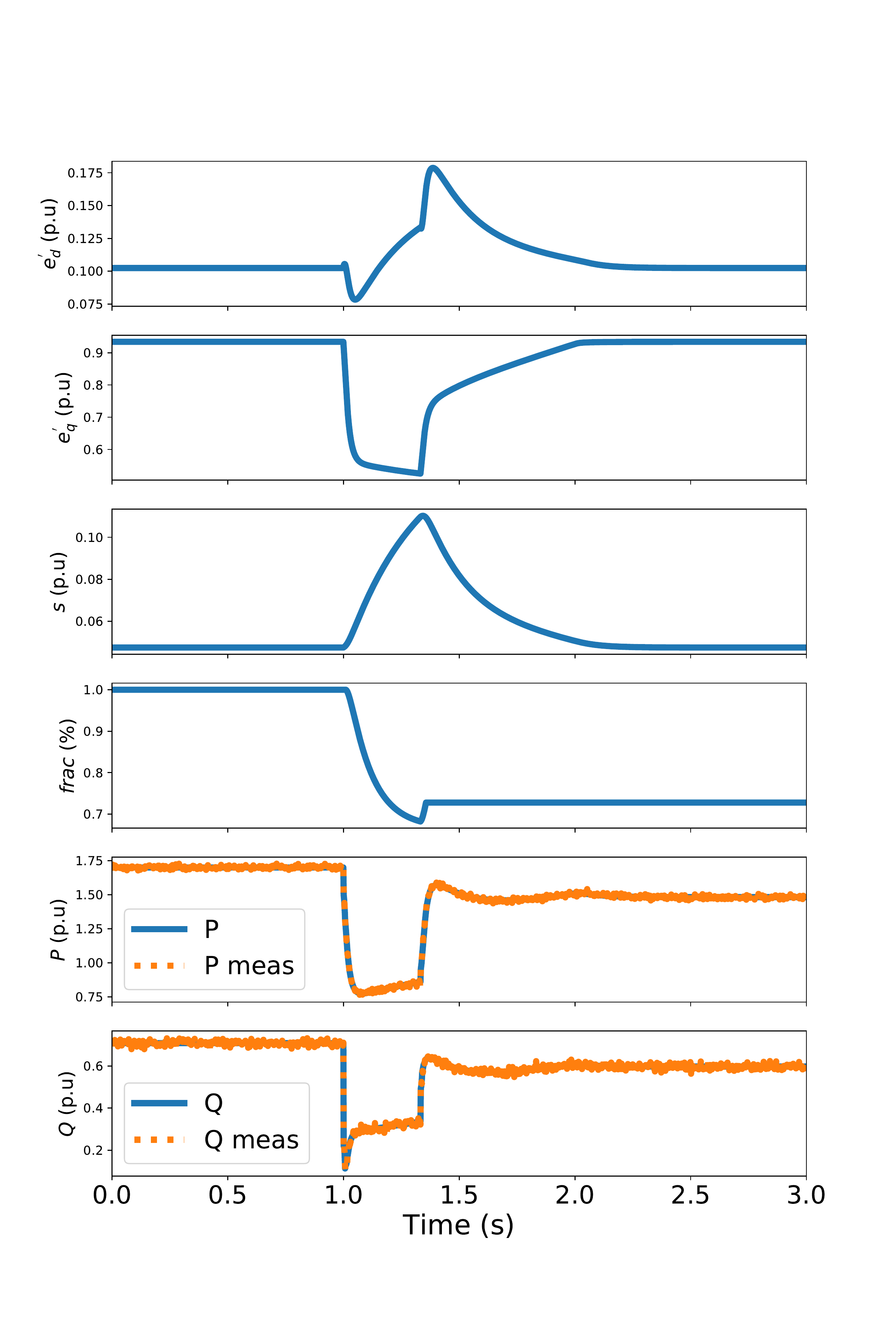}
\caption{States of the load when subjected to terminal voltage plotted in Fig.~\ref{fig:voltage_dist}. We measure active and reactive power at the terminals, where the noisy signal is plotted in discontinuous orange.}
\label{fig:load_states}
\end{figure}

Figures \ref{fig:voltage_dist} and \ref{fig:load_states} show the voltage signal and the resulting evolution of the states of the load and measurements. 

\section{Parameters}

Induction motor parameters: $r_a = 0.0138$ p.u, $x_a = 0.083$ p.u, $x_m = 3.0$ p.u, $r_1 = 0.055$ p.u, $x_1 = 0.053$ p.u, $H = 0.8$ sec., and initial value for initialization $P^0_{\textit{mot}} = 0.8$ p.u.

ZIP load parameters (in p.u.): $P_z = 0.6$, $P_i = 0.2$, $P_p = 0.1$, $Q_z = 0.2$, $Q_i = 0.05$, $Q_p = 0.05$.

Progressive tripping parameters: $v_{1\textit{off}} = 0.8$ p.u, $v_{2\textit{off}} = 0.2$ p.u, $Tr = 0.1$ s.

\end{appendices}

\bibliographystyle{ieeetr}
\bibliography{bibliography}

\end{document}